\documentclass[12pt]{amsart}
\usepackage{amsmath,amsthm,amsfonts,amssymb,times,hyperref}
\usepackage[top=1in, bottom=1in, left=1in, right=1in, footskip=0.5in]{geometry}
\usepackage{xcolor}
\usepackage{enumitem}

\theoremstyle{plain}
\newtheorem{theorem}{Theorem}[section]
\newtheorem{lemma}[theorem]{Lemma}

\theoremstyle{definition}
\newtheorem{comment}[theorem]{Comment}

\newtheorem{example}[theorem]{Example}

\DeclareMathOperator{\mmod}{\,mod\,}
\DeclareMathOperator{\lcm}{lcm}

\newcommand{\exactmid}{\mathrel\Vert}

\allowdisplaybreaks

\title{Minimum modulus problem for covering subsets of the integers}

\author{Michael Filaseta}
\address{Department of Mathematics, University of South Carolina, Columbia, SC 29208}
\email{filaseta@math.sc.edu}

\author{Alexandros Kalogirou}
\address{Erd{\H o}s Center, R\'enyi Institute, Budapest, Hungary H-1053}
\email{kalogira@email.sc.edu}

\author{Jonah Klein}
\address{Department of Mathematics, University of South Carolina, Columbia, SC 29208}
\email{jonah.klein@sc.edu}

\date{\today}


\begin{document}

\begin{abstract}
Expanding upon developments by B.~Hough and by P.~Balister, B.~Bollob{\'a}s, R.~Morris, J.~Sahasrabudhe and M.~Tiba, showing that the minimum modulus for a distinct covering system of the integers is bounded by $615999$, we establish that the minimum modulus for a distinct covering system for certain subsets of the integers is also bounded.  In particular, we establish that the minimum modulus for a distinct covering system of the prime numbers and the minimum modulus for a distinct covering system of the squares are bounded.  More generally, we show that if $f(x)$ is a non-constant polynomial in $\mathbb Z[x]$, then the minimum modulus for a distinct covering system of the set $\{ f(p): p \text{ prime} \}$ is bounded by a constant depending only on $f(x)$.
\end{abstract}

\maketitle

\section{Introduction}
A \textit{covering system} or \textit{covering} $\mathcal C$ for a set $\mathcal S \subseteq \mathbb Z$ is a finite set of residue classes $a \pmod{m}$ (or a finite set of congruences $x \equiv a \pmod{m}$), with each $m \in \mathbb Z$ exceeding $1$, such that every element of $\mathcal S$ belongs to at least one residue class in $\mathcal C$ (satisfies at least one congruence in $\mathcal C$).
The covering $\mathcal C$ of $\mathcal S \subseteq \mathbb Z$ is said to be a \textit{distinct} covering if the moduli $m$ of the residue classes (or congruences) in $\mathcal C$ are distinct.  

Work by B.~Hough \cite{hough} and later by P.~Balister, B.~Bollob{\'a}s, R.~Morris, J.~Sahasrabudhe and M.~Tiba \cite{bbmst3, bbmst2, bbmst} resolved a problem of P.~Erd{\H o}s \cite{erdos} and established that the minimum modulus in a distinct covering system of $\mathbb Z$ is bounded, where, due to \cite{bbmst}, we can take the bound to be $615999$.
A still open problem is to determine whether there exist distinct covering systems in which every modulus is an odd number greater than $1$.  Partial progress was made by B.~Hough and P.~Nielsen \cite{houghnielsen2019} who showed that no distinct covering system exists with every modulus exceeding $1$ and relatively prime to $6$.  P.~Balister, B.~Bollob{\'a}s, R.~Morris, J.~Sahasrabudhe and M.~Tiba~\cite{bbmst} improved on this result by showing that every distinct covering system, where each modulus exceeds $1$, has the least common multiple of the moduli divisible by either $2$, $9$ or $15$.  Inspired by such results, O.~Trifonov (see \cite{filhar}) asked what could be said about bounded minimal moduli or distinct coverings with odd moduli exceeding $1$ for distinct coverings of more general subsets $\mathcal S$ of $\mathbb Z$.

This topic was investigated in \cite{filhar}.  They showed, for example, that there is a distinct covering with odd moduli exceeding $1$ of the primes and of the numbers that are sums of two squares.  They also showed that the powers of $2$ and the Fibonacci numbers can both have distinct covering systems consisting of moduli that are only odd numbers exceeding any fixed bound.  The reader can consult \cite{filhar} for further results on this topic as well as the more recent paper \cite{bispelsetal}.

In this paper, we investigate sets $\mathcal S \subseteq \mathbb Z$ for which one can show that the minimum modulus of every distinct covering of $\mathcal S$ is bounded, where the bound depends on the set $\mathcal S$.  These results resolve, in particular, a number of questions posed in \cite{filhar}.  Our approach, based on ideas from \cite{bbmst}, leads to establishing that the minimum modulus in a distinct covering of the set of primes or in a distinct covering of the set of squares is bounded.  More generally, we obtain the following that includes such sets.

\begin{theorem}\label{corapplication3}
Let $f(x)$ be a non-constant polynomial in $\mathbb Z[x]$.  Then there is a constant $K = K(f)$ such that every distinct covering of $\{ f(p) : p \text{ prime} \}$ has a congruence with modulus $\le K$.
\end{theorem}

\section{Preliminaries and a general result}
We begin with a subset $\mathcal S$ of the integers which we would like to show has the property that any finite distinct covering $\mathcal C$ of $\mathcal S$ has a congruence with modulus $\le K$ where $K$ is some constant that depends only on $\mathcal S$ or more generally only on some property of $\mathcal S$.  
Recall that, in the previous section, we defined covering systems to be \textit{finite} sets of congruences.
For integers $a$ and $m$ with $m > 0$, denote by $a \mmod m$ the $b \in \{ 0, 1, \ldots, m-1 \}$ satisfying $b \equiv a \pmod{m}$.  
For each prime $p$ and each positive integer $\gamma$, we define
\begin{equation}\label{Apgammadef}
\mathcal A_{p^{\gamma}} = \{ s \mmod p^{\gamma} : s \in \mathcal S \}
\subseteq \{ 0, 1, \ldots, p^{\gamma}-1 \}.
\end{equation}
The sets $\mathcal A_{p^{\gamma}}$ are residue restrictions imposed by the elements of the set $\mathcal S$.  If there are no restrictions for a given prime $p$ and positive integer $\gamma$, that is no residue class is to be omitted modulo $p^{\gamma}$, then we will have $\mathcal A_{p^{\gamma}} = \{ 0, 1, \ldots, p^{\gamma}-1 \}$.  
Define
\[
\rho(p^{\gamma}) = \dfrac{|\mathcal A_{p^{\gamma}}|}{p^{\gamma}}.
\]
The following property and condition on the sets $\mathcal S$ are of interest to us.

\medskip \noindent
\textbf{Property A.}
\textit{For every $s \in \mathcal S$, for every prime $p$ and for every $\gamma \in \mathbb Z^{+}$, there is an $a \in \mathcal A_{p^{\gamma}}$ such that $s \equiv a \pmod{p^{\gamma}}$.}

\medskip \noindent
\noindent
\textbf{Condition A.}
\textit{There is a constant $c_{0} = c_{0}(\mathcal S) > 0$ such that $\rho(p^{\gamma}) \ge c_{0} $ for all primes $p$ and $\gamma \in \mathbb Z^{+}$.}

\medskip
That Property~A holds follows from the definition of $\mathcal A_{p^{\gamma}}$ in \eqref{Apgammadef}.  
In our arguments to come, we could relax Condition~A and instead require $\rho(p^{\gamma})$ to be larger than a suitable function of $p$ that tends to $0$ as $p \rightarrow \infty$.  In the end, we saw no advantage to making such a change, so we have kept Condition~A in the simpler form stated above.

\begin{example}
For $\mathcal S = \mathbb Z$, we have $\mathcal A_{p^{\gamma}} = \{ 0, 1, \ldots, p^{\gamma}-1 \}$ so that $\rho(p^{\gamma}) = 1$ for all primes $p$ and positive integers $\gamma$.  
\end{example}

\begin{example}
Fix $y \in \mathbb Z^{+}$.  
Let $f(x) \in \mathbb Z[x]$ of degree $n \ge 1$.  
Let $\mathcal S$ be the set of values $f(q)$, where $q$ varies over the primes greater than $y$.  
Then, for a prime $p$ and positive integer $\gamma$, we have
\[
\mathcal A_{p^{\gamma}} = \{ f(q) \mmod p^{\gamma} : q > y, q \text{ prime} \}.
\]  
We establish that there is a lower bound $c_{0} = c_{0}(f) > 0$ on $\rho(p^{\gamma})$.  The argument is reminiscent of arguments for Hensel lifting.

\begin{lemma}
Let $y \in \mathbb Z^{+}$, and let $f(x)$ be a non-constant polynomial in $\mathbb Z[x]$.  Then there is a constant $c_{0}(f) > 0$, independent of $y$, such that for every prime $p$ and every $\gamma \in \mathbb Z^{+}$, we have
\[
\big| \{ f(q) \mmod p^{\gamma} : q > y, q \text{ prime} \} \big| \ge c_{0}(f) \,p^{\gamma}.
\]
\end{lemma}

\begin{proof}
We use $H = H(f)$ to denote the height of $f(x)$, that is the maximum of the absolute values of the coefficients of $f(x)$.  
Let $n = \deg f \ge 1$.
We consider two cases depending on whether $p > \max\{2n,H\}$ or not.

First, fix a prime $p > \max\{2n,H\}$ and a positive integer $\gamma$.  
As $f'(x)$ is a non-zero polynomial modulo $p$ of degree $n-1$, it has at most $n-1$ distinct roots modulo $p$.  
For each integer $t$, we also have that $f(x)-t$ has at most $n$ roots modulo $p$.  
Since $p > 2n$, we deduce that, of the $\ge p - n$ different $m \in \{ 1, \ldots, p-1 \}$ for which $f'(m) \not\equiv 0 \pmod{p}$, there are 
\[
\ge \bigg\lceil \dfrac{p-n}{n} \bigg\rceil \ge \dfrac{p-n}{n} > \dfrac{p-(p/2)}{n} = \dfrac{p}{2n} 
\]
different values of $f(m)$ modulo $p$.  We extend each of these to $\ge p^{\gamma-1}$ different values of $f(m)$ modulo $p^{\gamma}$.  More precisely, let $m_{0}$ be such that $f'(m_{0}) \not\equiv 0 \pmod{p}$ and $m_{0} \not\equiv 0 \pmod{p}$, and suppose we know for some $j \ge 1$ that ($*$) holds where ($*$) is defined as follows.  

\begin{enumerate}[leftmargin=2em,rightmargin=2em]
\begin{item}[($*$)]
Let $S_{j}$ be the set of the $p^{j-1}$ integers $m \in \{ 0, 1, \ldots, p^{j}-1 \}$ that are congruent to $m_{0}$ modulo $p$.  Then the set $\{ f(m) \mmod p^{j} : m \in S_{j} \}$ also has size $p^{j-1}$.  
\end{item}
\end{enumerate}

\noindent
We show that ($*$) holds with $j$ replaced by $j+1$.  Then, since ($*$) holds with $j = 1$,
we can deduce that ($*$) holds for $j = \gamma$.

Fix $m'_{0} \in S_{j}$.  Then $f'(m'_{0}) \equiv f'(m_{0}) \not\equiv 0 \pmod{p}$.  
By the Maclaurin series for $f(x+m'_{0})$, we see that for every $k \in \mathbb Z$, we have
\[
f(m'_{0} + k p^{j}) \equiv f(m'_{0}) + kp^{j} f'(m'_{0}) \pmod{p^{j+1}}.
\]
Since $p \nmid f'(m'_{0})$, we see that for $k \in \{ 0, 1, \ldots, p-1 \}$, the values of $f(m'_{0}) + kp^{j} f'(m'_{0})$ are the same as $f(m'_{0})$ modulo $p^{j}$ and distinct modulo $p^{j+1}$ (the difference of any two of them is not divisible by $p^{j+1}$).  
Note that by the definition of $S_{j}$, we have $m'_{0} + k p^{j} \equiv m_{0} \pmod{p}$, so as the numbers $m'_{0} \in S_{j}$ and $k \in \{ 0, 1, \ldots, p-1 \}$ vary, the numbers $m'_{0} + k p^{j}$ are the $p^{j}$ elements of $S_{j+1}$.
Since the $p$ values of $f(m'_{0}) + kp^{j} f'(m'_{0})$ are congruent to $f(m'_{0})$ modulo $p^{j}$ and distinct modulo $p^{j+1}$ and the set $\{ f(m) \mmod p^{j} : m \in S_{j} \}$ has size $p^{j-1}$, we see that there are precisely $p^{j}$ elements in $\{ f(m) \mmod p^{j+1} : m \in S_{j+1} \}$.  This proves ($*$) by induction.  
Thus, each of the $\ge p/(2n)$ different values of $f(m_{0})$ modulo $p$ with $f'(m_{0}) \not\equiv 0 \pmod{p}$ and  $m_{0} \not\equiv 0 \pmod{p}$ extends to $p^{\gamma-1}$ different values of $f(m)$ modulo $p^{\gamma}$ which are the same as $f(m_{0})$ modulo $p$.  We deduce then that 
\begin{equation}\label{distvalsoff(x)eq1}
\big| \{ f(m) \mmod p^{\gamma} : m \in \mathbb Z, p \nmid m \} \big| \ge \dfrac{p^{\gamma}}{2n},
\quad \text{for all $p > \max\{2n,H\}$ and $\gamma \in \mathbb Z^{+}$}.
\end{equation}
To complete our analysis in the case that $p > \max\{2n,H\}$ and $\gamma \in \mathbb Z^{+}$, we observe that for each integer $m \not\equiv 0 \pmod{p}$, Dirichlet's Theorem (on primes in arithmetic progressions) implies that there is a prime $q > y$ satisfying $q \equiv m \pmod{p^{\gamma}}$ and, hence, $f(q) \equiv f(m) \pmod{p^{\gamma}}$.  Therefore, from \eqref{distvalsoff(x)eq1}, we see that
\begin{equation}\label{distvalsoff(x)eq1pt5}
\big| \{ f(q) \mmod p^{\gamma} : q > y, q \text{ prime} \} \big| \ge \dfrac{p^{\gamma}}{2n},
\quad \text{for all $p > \max\{2n,H\}$ and $\gamma \in \mathbb Z^{+}$}.
\end{equation}

Next, consider a prime $p \le \max\{2n,H\}$.  Since the number of such primes only depends on $f(x)$ and there is at least one value of $f(q)$ modulo $p^{\gamma}$ with $q >y$ prime, we can and do consider $\gamma$ to be large; more precisely, for $\gamma_{0} \in \mathbb Z^{+}$ fixed but arbitrarily large, we have
\begin{equation}\label{distvalsoff(x)eq2}
\begin{split}
\big| \{ f(q) \mmod p^{\gamma} : q > y, q &\text{ prime} \} \big| \ge 1 \ge \dfrac{p^{\gamma}}{\max\{2n,H\}^{\gamma_{0}}}, \\[6pt]
&\quad \text{for all primes $p \le \max\{2n,H\}$ and all $\gamma \in [1,\gamma_{0}] \cap \mathbb Z$}.
\end{split}
\end{equation}
We consider now $p \le \max\{2n,H\}$ and $\gamma > \gamma_{0}$, where $\gamma_{0}$ is large and to be determined.

Fix $m_{0} \in \mathbb Z$ satisfying $p \nmid m_{0}$ and $f'(m_{0}) \ne 0$, which is possible since $f'(x)$ has finitely many roots.  Let $s$ be the nonnegative integer for which $p^{s} \exactmid f'(m_{0})$.  Fix an integer $r > s$.  We take $\gamma_{0} = r+s$ above so that now we are interested in $\gamma > r+s$.  
Then an integer $m$ of the form $m_{0} + kp^{r}$ for $k \in \mathbb Z$ satisfies
\[
f'(m) \equiv f'(m_{0}) \pmod{p^{r}},
\]
which implies, since $r > s$, that $p^{s} \exactmid f'(m)$.  Also, since $2r \ge r+s+1$, we see that
\[
f(m) \equiv f(m_{0}) + kp^{r} f'(m_{0}) \pmod{p^{r+s+1}}.
\]
For $k \in \{ 0, 1, \ldots, p-1 \}$, we obtain $p$ values of $m = m_{0} + kp^{r}$ giving $p$ incongruent values of $f(m)$ modulo $p^{r+s+1}$.  Furthermore, these $p$ values of $m$ are congruent modulo $p^{r}$ but incongruent modulo $p^{r+1}$.  Replacing $m_{0}$ now with one of these $p$ values of $m$, say $m_{1}$, and replacing $r$ with $r+1$ in the above, we obtain $p$ values of $m = m_{1} + kp^{r+1}$ giving $p$ incongruent values of $f(m)$ modulo $p^{r+s+2}$ with $m \equiv m_{1} \pmod{p^{r+1}}$ and with the $p$ values of $m$ incongruent modulo $p^{r+2}$.  By an induction argument, somewhat similar to what was done for $p > \max\{2n,H\}$, we deduce that, for integers $\gamma \ge \gamma_{0}$, there are $\ge p^{\gamma - r - s}$ incongruent values of $f(m)$ modulo $p^{\gamma}$ with $m \equiv m_{0} \not\equiv 0 \pmod{p}$.  We deduce that
\begin{equation*}
\begin{split}
\big| \{ f(m) \mmod p^{\gamma} : m \in &\ \mathbb Z, p \nmid m \} \big| \ge p^{\gamma - r - s} = \dfrac{p^{\gamma}}{p^{r+s}} \ge \dfrac{p^{\gamma}}{\max\{2n,H\}^{r+s}}\\[6pt]
&\quad \text{for all primes $p \le \max\{2n,H\}$ and all $\gamma \in (\gamma_{0},\infty) \cap \mathbb Z$}.
\end{split}
\end{equation*}
As before we deduce then, as a consequence of Dirichlet's Theorem, that
\begin{equation}\label{distvalsoff(x)eq3}
\begin{split}
\big| \{ f(q) \mmod p^{\gamma} : q > y, &\ q \text{ prime} \} \big| \ge \dfrac{p^{\gamma}}{\max\{2n,H\}^{r+s}}\\[6pt]
&\quad \text{for all primes $p \le \max\{2n,H\}$ and all $\gamma \in (\gamma_{0},\infty) \cap \mathbb Z$}.
\end{split}
\end{equation}
Combining \eqref{distvalsoff(x)eq1pt5}, \eqref{distvalsoff(x)eq2} and \eqref{distvalsoff(x)eq3}, the lemma follows. 
\end{proof}

For $y \ge 1$, we refer to a $y$-smooth number as a positive integer having each of its prime factors $\le y$ and to a $y$-rough number as a positive integer having a minimal prime factor that exceeds $y$.  We will examine coverings for both of these sets. 
We note that the primes $q > y$ are $y$-rough numbers, so a simple corollary of the lemma above that will play a role later in the paper is the following.

\begin{lemma}\label{distvalsoff(x)}
Let $y \in \mathbb Z^{+}$, and let $f(x)$ be a non-constant polynomial in $\mathbb Z[x]$.  Then there is a constant $c_{0}(f) > 0$, independent of $y$, such that for every prime $p$ and every $\gamma \in \mathbb Z^{+}$, we have
\[
\big| \{ f(m) \mmod p^{\gamma} : m \in \mathbb Z^{+}, m \text{ is $y$-rough} \} \big| \ge c_{0}(f) \,p^{\gamma}.
\]
\end{lemma}
\end{example}

For Theorem~\ref{corapplication3}, we seemingly want to take $\mathcal S = \{ f(p) : p \text{ prime} \}$.  However, we did not see a way to directly approach such a choice for $\mathcal S$.  
Property~A and Condition~A on the set $\mathcal S$ are not going to be sufficient for our arguments to show that the minimum modulus is bounded, and an additional condition we want to impose does not allow us to apply our general result to the set $\{ f(p) : p \text{ prime} \}$ directly.  
We will begin by focusing on the choice
\begin{equation}\label{Sforf(m)rough}
\mathcal S = \{ f(m) : m \in \mathbb Z^{+}, m \text{ is $y$-rough} \},
\end{equation}
for a fixed $y \in \mathbb Z^{+}$.  We are initially interested then in showing that the minimum modulus is bounded for congruences in a distinct covering of $\mathcal S$ as in \eqref{Sforf(m)rough}.  

Observe that, despite Property~A, it is not a consequence of \eqref{Apgammadef} that if $s \in \mathbb Z^{+}$ such that for every prime $p$ and every positive integer $\gamma$ we have $s \equiv a \pmod{p^{\gamma}}$ for some $a \in \mathcal A_{p^{\gamma}}$, then $s \in \mathcal S$.  By way of a simple example, consider the set 
\begin{equation}\label{sprimeeg}
\mathcal S' = \{ n \in \mathbb Z^{+} : n = p^{\gamma} + 5 \text{ for some prime $p$ and some $\gamma \in \mathbb Z^{+}$} \}.
\end{equation} 
By \eqref{Apgammadef}, independent of $p$ and $\gamma$, we see that $5 \equiv a \pmod{p^{\gamma}}$ for some $a \in A_{p^{\gamma}}$, but $5 \not\in \mathcal S'$.   

Our goal will be to show that the following additional condition on $\mathcal S$ is sufficient for showing the minimum modulus for congruences in a distinct covering of $\mathcal S$ is bounded.

\medskip
\noindent
\textbf{Condition B.}
\textit{For $j$ a positive integer, let $p_{j}$ denote the $j^{\rm th}$ prime.  
For $r$ a positive integer and 
\[
\vec{v} = \langle \gamma_{1}, \ldots, \gamma_{r} \rangle,
\]
where the $\gamma_{j}$ are positive integers, set 
\[
\mathcal S_{r}(\vec{v}) = \big\{ n \in \mathbb Z : \text{for every } j \in [1,r] \cap \mathbb Z, \text{we have } n \equiv a {\hskip -8pt}\pmod{p_{j}^{\gamma_{j}}} \text{ for some } a \in \mathcal A_{p_{j}^{\gamma_{j}}} \big\}.
\]
There is a constant $c_{1} > 0$ such that, 
for every $r$ and $\vec{v}$ as above with $\min\{ \gamma_{1}, \ldots, \gamma_{r} \} \ge c_{1}$,  and for every $n \in \mathcal S_{r}(\vec{v})$, there is an $n' \in \mathcal S$ such that $n' \equiv n \pmod{p_{j}^{\gamma_{j}}}$ for every integer $j \in [1,r]$.}

\medskip
Before continuing, we justify that there is some non-constant $f(x) \in \mathbb Z[x]$ such that $\mathcal S =  \{ f(p) : p \text{ prime} \}$ does not satisfy Condition~B, but for all non-constant $f(x) \in \mathbb Z[x]$ and positive integers $y$, the set $\mathcal S$ in \eqref{Sforf(m)rough} does satisfy Condition~B.

First, suppose $\mathcal S =  \{ f(p) : p \text{ prime} \}$ where $f(x) = x$.  
Fix integers $r \ge 2$ and $c_{1} > 0$.  Let $\gamma = \max\{ 2, \lceil c_{1} \rceil \}$.  
Take $\gamma_{j} = \gamma$ for all integers $j \in [1,r]$.
Let $n$ be a positive integer satisfying the congruences $x \equiv 2 \pmod{2^{\gamma}}$, 
$x \equiv 3 \pmod{3^{\gamma}}$, and $x \equiv 1 \pmod{p_{j}^{\gamma}}$ for $3 \le j \le r$.
By Dirichlet's Theorem in the case $j \ge 3$, we see that for each fixed $j \in [1,r]$, there is at least one prime satisfying the congruence modulo $p_{j}^{\gamma}$ just described, so it follows that $n \equiv a \pmod{p_{j}^{\gamma}}$ for some $a \in \mathcal A_{p_{j}^{\gamma}}$.  Hence, $n \in \mathcal S_{r}(\langle \gamma, \ldots, \gamma \rangle)$.  On the other hand, if $n' \equiv n \pmod{p_{j}^{\gamma}}$ for every $j \in [1,r]$, then $n' \equiv 2 \pmod{2^{\gamma}}$ and $n' \equiv 3 \pmod{3^{\gamma}}$ implies that $6$ divides $n'$.  Hence, there is no element $n'$ of $\mathcal S =  \{ f(p) : p \text{ prime} \} =  \{ p : p \text{ prime} \}$ satisfying $n' \equiv n \pmod{p_{j}^{\gamma}}$ for every integer $j \in [1,r]$.  Therefore, Condition~B does not hold if $\mathcal S =  \{ p : p \text{ prime} \}$.  

Now, let $f(x)$ be a non-constant polynomial in $\mathbb Z[x]$ and $y$ be a positive integer.  Let $\mathcal S$ be as in \eqref{Sforf(m)rough}.  We take $c_{1} = 1$.  Fix a positive integer $r$ and a vector $\vec{v} = \langle \gamma_{1}, \ldots, \gamma_{r} \rangle$ as in Condition~B.  Let $n \in \mathcal S_{r}(\vec{v})$.  
Then, for each $j \in \{ 1, 2, \ldots, r \}$, there is a $y$-rough number $m_{j}$ such that $f(m_{j}) \equiv n \pmod{p_{j}^{\gamma_{j}}}$.  If there is a $p_{j}$ dividing $m_{j}$ for some $j \in \{ 1, 2, \ldots, r \}$, then let $j_{0}$ denote the minimal such $j$.  Otherwise, set $j_{0} = r+1$.  Observe that if $j_{0} \le r$, then since $m_{j_{0}} \equiv 0 \pmod{p_{j_{0}}}$ and $m_{j_{0}}$ is $y$-rough, we have $p_{j_{0}} > y$. Thus, $p_{j} > y$ for all $j \ge j_{0}$.  
Using $\nu_{p_{j}}$ as the usual $p_{j}$-adic valuation, we define
\[
P_{0} = \prod_{j = 1}^{r} p_{j}^{\nu_{p_{j}}(m_{j})} = \prod_{j = j_{0}}^{r} p_{j}^{\nu_{p_{j}}(m_{j})}
\qquad \text{and} \qquad
m_{0} = p' P_{0},
\]
where $p' > \max\{ p_{r}, y \}$ is a prime satisfying
\[
p' \cdot \dfrac{P_{0}}{p_{j}^{\nu_{p_{j}}(m_{j})}} \equiv \dfrac{m_{j}}{p_{j}^{\nu_{p_{j}}(m_{j})}} \pmod{p_{j}^{\gamma_{j}}},
\qquad \text{for every $j \in \{ 1, 2, \ldots, r \}$},
\]
which exists by the Chinese Remainder Theorem and Dirichlet's Theorem.
Since all the prime divisors of $P_{0}$ and $p'$ are greater than $y$, we see that $m_{0}$ is $y$-rough.  
Also, $m_{0} = p' P_{0} \equiv m_{j} \pmod{p_{j}^{\gamma_{j}}}$ for every $j \in \{ 1, 2, \ldots, r \}$.   
Thus, $f(m_{0}) \equiv f(m_{j}) \equiv n \pmod{p_{j}^{\gamma_{j}}}$ for every $j \in \{ 1, 2, \ldots, r \}$.  
We deduce then that $n' = f(m_{0}) \in \mathcal S$ and $n' \equiv n \pmod{p_{j}^{\gamma_{j}}}$ for every $j \in \{ 1, 2, \ldots, r \}$, showing that Condition~B holds for $\mathcal S$ as in \eqref{Sforf(m)rough}.

\begin{comment}
With $\mathcal S = \mathcal S'$ as in \eqref{sprimeeg}, if $r \ge 3$ in Condition~B, then $5 \in S_{r}(\vec{v})$ but one can show as follows that there is no $n' \in \mathcal S$ such that $n' \equiv 5 \pmod{p_{j}^{\gamma_{j}}}$ for every integer $j \in [1,r]$.  Assume such an $n' \in \mathcal S$ exists.  Since $r \ge 3$, we have $n' \equiv 5 \pmod{5}$ so that $n'$ is divisible by $5$.  However, the only elements of $\mathcal S$, that is the only numbers of the form $p^{\gamma}+5$, that are divisible by $5$ are the numbers $5^{\gamma} + 5$ for some positive integer $\gamma$.  Since $5^{\gamma} + 5 \not\equiv 5 \pmod{2}$, we obtain a contradiction to $n' \equiv 5 \pmod{p_{1}^{\gamma_{1}}}$.  We deduce that $\mathcal S$ does not satisfy Condition~B in this case.
\end{comment}

We are ready now to state our main result from which Theorem~\ref{corapplication3} will follow.

\begin{theorem}\label{densitythm}
Given the above setup with $\mathcal S$ satisfying Conditions~A and B, there is a constant $K = K(c_{0})$ such that every distinct covering $\mathcal C$ of $\mathcal S$ has a congruence with modulus $\le K$.  
\end{theorem}

Before addressing the proof of Theorem~\ref{densitythm}, we note that the lemma below follows from Theorem~\ref{densitythm} and our prior discussion.

\begin{lemma}\label{corapplication2}
Let $y \in \mathbb Z^{+}$, and let $f(x)$ be a non-constant polynomial in $\mathbb Z[x]$.  
There is a $K = K(f)$, independent of $y$, such that every distinct covering of the set $\mathcal S$ in \eqref{Sforf(m)rough} has a congruence with modulus $\le K$.
\end{lemma}

With a little more work, we also obtain a proof of Theorem~\ref{corapplication3}.

\begin{proof}[Proof of Theorem~\ref{corapplication3}]
Let $\mathcal C$ be a distinct covering of the set $\{ f(p) : p \text{ prime} \}$.  
Let $r$ and $\gamma_{j}$ be positive integers such that the moduli appearing in congruences in $\mathcal C$ are all divisors of $L = p_{1}^{\gamma_{1}} \cdots p_{r}^{\gamma_{r}}$, where as before $p_{j}$ denotes the $j^{\rm th}$ prime.  Take $y = p_{r}$.  To establish Theorem~\ref{corapplication3}, it suffices by Lemma~\ref{corapplication2} to show that $\mathcal C$ is a covering of the set $\mathcal S$ in \eqref{Sforf(m)rough}.
Let $m'$ be a $y$-rough number.  
Since $m'$ is not divisible by a prime $\le y = p_{r}$, we see that $m'$ is relatively prime to $L$.  By Dirichlet's Theorem, there is a prime $q \equiv m' \pmod{L}$.  
As $\mathcal C$ is a covering of the set $\{ f(p) : p \text{ prime} \}$, we obtain that $f(q)$ satisfies some congruence, say $x \equiv a \pmod{m}$, in $\mathcal C$.  Now, $q \equiv m' \pmod{L}$ and $m \mid L$ implies $f(m') \equiv f(q) \equiv a \pmod{m}$, from which we deduce $f(m')$ satisfies a congruence in $\mathcal C$.  The theorem follows.
\end{proof}

We elaborate on consequential examples of sets $\mathcal S \subseteq \mathbb Z$ where distinct coverings have the minimal modulus among their congruences bounded.
As a distinct covering of the squares is also a distinct covering of the squares of the primes, we see that
Theorem~\ref{corapplication3} implies that there is a $K$ such that every distinct covering of the squares includes a congruence with modulus $\le K$.  This is related to a question posed in \cite{filhar} where it was shown that there is an odd covering (a distinct covering where all moduli are odd and $> 1$) for numbers that are sums of two squares.  The question was whether the minimum modulus for congruences in a distinct covering of the set of integers that are a sum of two squares is bounded above.  Since the set of the sum of two squares includes the set of squares, we see that such a covering will in fact have a bounded minimum modulus.  For an integer $k \ge 2$, Theorem~\ref{corapplication3} also implies that distinct coverings of the set of $k$-full numbers (since they include the set of $k^{\rm th}$ powers of primes) have a bounded minimum modulus.  Also, Theorem~\ref{corapplication3} with $f(x) = x$ implies that distinct coverings of the set of squarefree numbers (and consequently the set of $k$-free numbers for an integer $k \ge 2$) have a bounded minimum modulus.  As another example of an application of Theorem~\ref{corapplication3}, we note that, despite the fact that $\mathcal S = \mathcal S'$ defined by \eqref{sprimeeg} does not satisfy Condition~B, the minimal modulus of the congruences in a distinct covering of $\mathcal S'$ is bounded.  To see this, observe that if $\mathcal C$ is a distinct covering of $\mathcal S'$, then the covering $\mathcal C'$ obtained by replacing each congruence $x \equiv a \pmod{m}$ in $\mathcal C$ by $x \equiv a-5 \pmod{m}$ is a distinct covering of the primes.  Thus, by Theorem~\ref{corapplication3}, the minimal modulus of $\mathcal C'$ and, hence, $\mathcal C$ is bounded.

As noted in the introduction, in the opposite direction, from \cite{filhar}, one can show that the minimum modulus can be arbitrarily large for covering systems of the powers of $2$.  This result generalizes to obtain the following.

\begin{theorem}\label{ysmooththm}
Fix $y \in \mathbb Z^{+}$.  Then for every positive integer $M$, there is a distinct covering of the $y$-smooth numbers with the minimum modulus greater than $M$. 
\end{theorem}

\begin{proof}
Let $p_1, \ldots, p_r$ be the primes $\le y$.  Let $e_1, \ldots, e_r$ be positive integers for which $p_{j}^{e_{j}} > M$ for each $j \in \{1, 2, \ldots, r\}$.  We form a distinct covering of the $y$-smooth numbers by first including the residue classes $0 \pmod{p_{j}^{e_{j}}}$ in $\mathcal C$ for $j \in \{1, 2, \ldots, r\}$.  Every $y$-smooth number belongs to one of these residue classes except the $s = e_{1} \cdots e_{r}$ numbers of the form 
\begin{equation}\label{smootheq}
p_{1}^{u_{1}} \cdots p_{r}^{u_{r}}
\qquad \text{with } 0 \le u_{j} < e_{j}
\text{ for } j \in \{1, 2, \ldots, r\}.
\end{equation}
We choose arbitrary distinct positive integers $q_{1}, \ldots, q_{s}$ each greater than $M$ and different from $p_{1}^{e_{1}}, \ldots, p_{r}^{e_{r}}$. For each positive integer $a$ as in \eqref{smootheq}, we include in $\mathcal C$ a congruence class $a \pmod{q}$ where $q$ is chosen to be in $\{q_1,\ldots,q_s\}$ and different for each $a$.  Such a $\mathcal C$ is a distinct covering of the $y$-smooth numbers with the minimum modulus greater than $M$.
\end{proof}

There is a lot of freedom in the above proof.  For example, one can take the $e_{j}$ and then the $q_{j}$ to be arbitrarily large so that one can also deduce that a distinct covering of the $y$-smooth numbers exists where the sum of the reciprocals of the moduli in the covering is less than $\varepsilon$ for any prescribed $\varepsilon > 0$.

\section{The proof of Theorem~\ref{densitythm}}
Our next lemma switches the problem from one with infinitely many congruence restrictions on the elements of $\mathcal S$ given by Property~A (one for every prime $p$ and $\gamma \in \mathbb Z^{+}$) to one with finitely many such restrictions corresponding to covering a set $\mathcal S_{r}(\vec{v})$ as defined in Condition~B.  This is where Condition~B comes into play in the argument for Theorem~\ref{densitythm}.  Observe that \eqref{Apgammadef} implies that $\mathcal S \subseteq \mathcal S_{r}(\vec{v})$, so a covering of $S_{r}(\vec{v})$, for any $\vec{v}$, is necessarily a covering of $\mathcal S$.  The proof of the next lemma is showing that, under Condition~B, the opposite is also true, that a covering of $\mathcal S$ is also a covering of $\mathcal S_{r}(\vec{v})$.  We note that the argument is similar to the proof given for Theorem~\ref{corapplication3}.

\begin{lemma}\label{Srsuffices}
Let $\mathcal S$ satisfy Conditions~A and B.  
For a positive integer $r$ and a vector $\vec{v} = \langle \gamma_{1}, \ldots, \gamma_{r} \rangle$ with each $\gamma_{j} \in \mathbb Z^{+}$, set
\[
L_{r} 
= L_{r}(\vec{v})
= p_{1}^{\gamma_{1}} p_{2}^{\gamma_{2}} \cdots p_{r}^{\gamma_{r}}.
\]  
Suppose that there is a $K$, depending on $\mathcal S$, such that for every positive integer $r$ and every $\vec{v} = \langle \gamma_{1}, \ldots, \gamma_{r} \rangle$ with each $\gamma_{j} \in \mathbb Z^{+}$,  we have that every distinct covering $\mathcal C$ of the set $\mathcal S_{r}(\vec{v})$ where each modulus of a congruence in $\mathcal C$ divides $L_{r}$ 
is such that some modulus of a congruence in $\mathcal C$ is $\le K$.   Then for every distinct covering $\mathcal C$ of the set $\mathcal S$, there is a modulus of a congruence in $\mathcal C$ that is $\le K$.
\end{lemma}

\begin{proof}
Let $\mathcal C$ be a distinct covering of $\mathcal S$.  
Then there is a positive integer $r$ and a $\vec{v} = \langle \gamma_{1}, \ldots, \gamma_{r} \rangle$ with each $\gamma_{j} \in \mathbb Z^{+}$ such that every modulus of a congruence in $\mathcal C$ divides $L_{r} = L_{r}(\vec{v})$.  
We suppose as we can that each $\gamma_{j} \ge c_{1}$, where $c_{1}$ is as in Condition~B.
It suffices to show that $\mathcal C$ is also a covering of $\mathcal S_{r}(\vec{v})$.  
We make use of Condition~B.  If $n \in \mathcal S_{r}(\vec{v})$, then there is an $n' \in \mathcal S$ such that $n' \equiv n \pmod{p_{j}^{\gamma_{j}}}$ for every integer $j \in [1,r]$.  
In other words, we have $n' \equiv n \pmod{L_{r}}$.
As $\mathcal C$ is a covering of $\mathcal S$, we see that $n'$ satisfies a congruence in $\mathcal C$.  Since $n' \equiv n \pmod{L_{r}}$ and since every modulus of $\mathcal C$ divides $L_{r}$, we see that $n$ satisfies the same congruence in $\mathcal C$.  Since this is true for every $n \in \mathcal S_{r}(\vec{v})$, we deduce that $\mathcal C$ is a covering of $\mathcal S_{r}(\vec{v})$ as desired.
\end{proof}

We describe next a variant of the distortion method having origins in work of  B.~Hough~\cite{hough} and more fully developed by P.~Balister, B.~Bollob{\'a}s, R.~Morris, J.~Sahasrabudhe and M.~Tiba~\cite{bbmst3,bbmst2,bbmst}.  With $r$ and $\vec{v} = \langle \gamma_{1}, \ldots, \gamma_{r} \rangle$ fixed as defined in  Lemma~\ref{Srsuffices}, our goal is to show that every distinct covering $\mathcal C$ of $\mathcal S_{r}(\vec{v})$ with moduli dividing $L_{r} = L_{r}(\vec{v})$ has a modulus that is $\le K$, where $K$ depends only on $c_{0}$ in Condition~A.  One of the basic ideas in the distortion method is the assignment of a probability distribution, or weights, to the set $[1,L_{r}] \cap \mathcal S_{r}(\vec{v})$.  The goal is to show that if $\mathcal C'$ is a set of congruences with distinct moduli dividing $L_{r}$ and the minimum modulus of $\mathcal C'$ is too large, then the sum of the weights of the integers in $[1,L_{r}] \cap \mathcal S_{r}(\vec{v})$ that satisfy no congruence in $\mathcal C'$ is positive.  Thus, if the minimum modulus is too large, then $\mathcal C'$ cannot be a covering of $\mathcal S_{r}(\vec{v})$.  To simplify notation below, we set
\[
[ a, m ] = \{ n \in \mathbb Z : n \equiv a {\hskip -8pt}\pmod{m} \},
\]
so we can use $[a,m]$ when referring to the elements of a congruence class covered by $\mathcal C'$.  With some abuse of notation, we will write $[a,m] \in \mathcal C'$ to indicate that the congruence $x \equiv a \pmod{m}$ is in $\mathcal C'$.  
We also make use of the notation $P(m)$ to denote the largest prime factor of an integer $m > 1$ and $P(1) = 1$.

We suppose throughout, as we may, that each congruence in $\mathcal C'$ is satisfied by some element of $\mathcal S_{r}(\vec{v})$; otherwise, we could simply remove the congruence from $\mathcal C'$ without affecting whether $\mathcal C'$ is a covering of $\mathcal S_{r}(\vec{v})$ and the conclusion of Theorem~\ref{densitythm}.

By the Chinese Remainder Theorem, there is a one-to-one correspondence between the integers in $[1,L_{r}] \cap \mathcal S_{r}(\vec{v})$ and their residue classes modulo $p_{i}^{\gamma_{i}}$ which we write as
\[
[1,L_{r}] \cap \mathcal S_{r}(\vec{v}) \longleftrightarrow \mathcal A_{p_{1}^{\gamma_{1}}} \times \mathcal A_{p_{2}^{\gamma_{2}}} \times \cdots \times \mathcal A_{p_{r}^{\gamma_{r}}}.
\]
We set $Q_{0} = \{ 1 \}$ and
\[
Q_{i} = \mathcal A_{p_{1}^{\gamma_{1}}} \times \mathcal A_{p_{2}^{\gamma_{2}}} \times \cdots \times \mathcal A_{p_{i}^{\gamma_{i}}}, \qquad \text{for $1 \le i \le r$},
\]
and view elements of $Q_{i}$, for $i \ge 1$, as $(x,y)$ where $x \in Q_{i-1}$ and $y \in \mathcal A_{p_{i}^{\gamma_{i}}}$.  In particular, in the case $i = 1$ where $x \in Q_{0} = \{ 1\}$, the tuple $(x,y) = (1,y)$ is equated with $y \in \mathcal A_{p_{1}^{\gamma_{1}}}$.
We can further view $Q_{i}$ as representing a subset of the integers in $[1,\ell_{i}]$, where
$\ell_{i} = p_{1}^{\gamma_{1}} \cdots p_{i}^{\gamma_{i}}$.  An integer $a$ is in this subset of $[1,\ell_{i}]$ precisely when there is at least one integer $b$ such that if $n \equiv a \pmod{\ell_{i}}$ and $n \equiv b \pmod{\ell_{r}/\ell_{i}}$, then $n \in \mathcal S_{r}(\vec{v})$.

We let $B_{i}$ denote the integers in $Q_{i}$ which satisfy at least one congruence in $\mathcal C'$ where the modulus of the congruence is divisible by $p_{i}$ and $p_{i}$-smooth.  From our notation above, we can write
\[
B_{i} = \{ n \in Q_{i} : n \in [a,m] \text{ for some } [a,m] \in \mathcal C' \text{ where } P(m) = p_{i} \}.
\]
Similar to \cite{bbmst}, we define a probability distribution $\mathbb P_{i}$ on the elements of $Q_{i}$, depending on a parameter $\delta_{i} \in [0,1/2]$ to be chosen later.  We set $\mathbb P_{0}(1) = 1$.  
For $i \ge 1$ and $x \in Q_{i-1}$, we define
\[
\alpha_i(x) 
= \dfrac{|\{y \in \mathcal A_{p_{i}^{\gamma_{i}}} : (x,y)\in B_i\}|}{|\mathcal A_{p_{i}^{\gamma_{i}}}|},
\]
and, for fixed $y \in \mathcal A_{p_{i}^{\gamma_{i}}}$, define
\begin{equation}\label{probdistPi}
\mathbb{P}_i(x,y)=
\begin{cases} 
\max\bigg\{ 0,\dfrac{\alpha_i(x)-\delta_i}{\alpha_i(x)(1-\delta_i)}   \bigg\}\cdot \dfrac{\mathbb{P}_{i-1}(x)}{|\mathcal A_{p_{i}^{\gamma_{i}}}|}, & \text{if }(x,y)\in B_i, \\[15pt]
\min \bigg\{ \dfrac{1}{1-\alpha_i(x)},\dfrac{1}{1-\delta_i} \bigg\}\cdot \dfrac{\mathbb{P}_{i-1}(x)}{|\mathcal A_{p_{i}^{\gamma_{i}}}|}, & \text{if }(x,y)\notin B_i. 
\end{cases}
\end{equation}
These probability distributions $\mathbb{P}_i$ on $Q_{i}$ can be extended to probability distributions $\mathbb{P}^{\uparrow}_i$ on $Q_{r}$ (or $[1,L_{r}])$ uniformly so that for $x \in Q_{i}$ and $y \in \mathcal A_{p_{i+1}^{\gamma_{i+1}}} \times \cdots \times \mathcal A_{p_{r}^{\gamma_{r}}}$, we have $\mathbb P^{\uparrow}_{i}(x,y) = P_{i}(x)/(|\mathcal A_{p_{i+1}^{\gamma_{i+1}}}| \cdots |\mathcal A_{p_{r}^{\gamma_{r}}}|)$.  
We will take an initial string of $\delta_{i}$ to be $0$, so we note here that if $\delta_{1} = \delta_{2} = \cdots = \delta_{i_{0}} = 0$, then for each $i \in \{1, 2, \ldots, i_{0} \}$, we have $\mathbb P^{\uparrow}_{i}$ as defined above is the uniform distribution on $Q_{r}$; in other words, for such $i$ and every integer $x \in Q_{r}$, we have 
\[
\mathbb P^{\uparrow}_{i}(x) = \dfrac{1}{|Q_{r}|} = \dfrac{1}{|\mathcal A_{p_{1}^{\gamma_{1}}}| \,|\mathcal A_{p_{2}^{\gamma_{2}}}| \cdots |\mathcal A_{p_{r}^{\gamma_{r}}}|}.
\] 
We are interested in $\mathbb P = \mathbb P_{r}$, which one can check satisfies
\begin{equation}
\label{PQrequation}
\mathbb P(Q_{r}) = \sum_{x \in Q_{r}} \mathbb P(x) = \sum_{x \in Q_{r}} \mathbb P^{\uparrow}_{i}(x) = \sum_{x \in Q_{i}} \mathbb P_{i}(x) = \mathbb P^{\uparrow}_{i}(Q_{r}) = 1 
\qquad \text{for $1 \le i \le r$.}
\end{equation}
For a fixed $x \in Q_{i-1}$ with $i \ge 1$, we also note that
\begin{equation}
\label{PQiQiminus1eq}
 \sum_{y \in \mathcal A_{p_{i}^{\gamma_{i}}}} \mathbb P_{i}(x,y) = \mathbb P_{i-1}(x).
\end{equation}
From \eqref{PQiQiminus1eq}, one can see that the probability measure assigned by $\mathbb P$ to the set of all integers $n \in Q_{r}$ such that $n$ satisfies some congruence in $\mathcal C'$ with modulus $m$ where $P(m) = p_{i}$ is the same as $\mathbb P_{i}(B_{i})$.  Thus,
\begin{equation}\label{ultimategoal}
\mathbb P\big(\{ n \in Q_{r} : n \in [a,m] \text{ for some } [a,m] \in \mathcal C' \}\big) \le \sum_{i=1}^{r} \mathbb P_{i}(B_{i}).
\end{equation}

Observe that
\begin{align*}
\mathbb{P}_i(B_i)
&= \sum_{x \in Q_{i-1}} \sum_{\substack{y \in \mathcal A_{p_{i}^{\gamma_{i}}} \\ (x,y) \in B_{i}}} \mathbb{P}_i(x,y) \\
&= \sum_{x \in Q_{i-1}} \text{max}\bigg\{ 0,\dfrac{\alpha_i(x)-\delta_i}{\alpha_i(x)(1-\delta_i)}\bigg\}
\cdot \dfrac{\mathbb{P}_{i-1}(x)}{|\mathcal A_{p_{i}^{\gamma_{i}}}|} \cdot \sum_{\substack{y \in \mathcal A_{p_{i}^{\gamma_{i}}} \\ (x,y) \in B_{i}}} 1 \\
&=\dfrac{1}{1-\delta_i}\sum_{x \in Q_{i-1}} \text{max}\Big\{0,\alpha_i(x)-\delta_i\Big\}\cdot \mathbb{P}_{i-1}(x).
\end{align*}
For $\delta_{i} > 0$, the inequality 
\[
\max\{ 0, \alpha_i(x)-\delta_i \} \le \dfrac{\alpha^{2}_i(x)}{4\delta_{i}}
\] 
follows from $t^2-4t+4 = (t - 2)^{2} \ge 0$ with $t = \alpha_{i}(x)/\delta_{i}$.  
For $f: Q_{i-1} \to \mathbb R$, set
\[
\mathbb{E}_{i-1}[f] = \sum_{x \in Q_{i-1}} f(x)\,\mathbb{P}_{i-1}(x).
\]
Then we obtain
\begin{equation}\label{measureBi}
\mathbb{P}_i(B_i) \le \dfrac{1}{4\,\delta_i (1-\delta_i)} \sum_{x \in Q_{i-1}} \alpha_i^2(x)\,\mathbb{P}_{i-1}(x) =  \dfrac{1}{4\,\delta_i (1-\delta_i)} \,\mathbb{E}_{i-1}[\alpha_i^2].
\end{equation}
We want to adjust the bound for $\mathbb{E}_{i-1}[\alpha_i^2]$ given in \cite[Lemmas 3.6-3.7]{bbmst} for our current situation.  The adjustments are minor, but we will spell out the details for clarification.  
Observe that
\[
\alpha_{i}(x) = \dfrac{1}{|\mathcal A_{p_{i}^{\gamma_{i}}}|} \sum_{\substack{y \in \mathcal A_{p_{i}^{\gamma_{i}}} \\ (x,y) \in B_{i}}} 1
\le  \dfrac{1}{|\mathcal A_{p_{i}^{\gamma_{i}}}|} \sum_{y \in \mathcal A_{p_{i}^{\gamma_{i}}}}
\sum_{\substack{[a, m] \in \mathcal C' \\ P(m) = p_{i} \\ (x,y) \in [a,m]}} 1
=  \dfrac{1}{|\mathcal A_{p_{i}^{\gamma_{i}}}|} \sum_{\substack{[a, m] \in \mathcal C' \\ P(m) = p_{i}}}
\sum_{\substack{y \in \mathcal A_{p_{i}^{\gamma_{i}}} \\ (x,y) \in [a,m]}} 1.
\]
From the notation described above, the meaning here of $(x,y) \in [a,m]$ is that the integer in $[1,\ell_{i}]$ associated with $(x,y)$ belongs to the congruence class $a \pmod{m}$.  As $P(m) = p_{i}$, we may write $m = p_{i}^{j} m'$ where $P(m') < p_{i}$ and $1 \le j \le \gamma_{i}$.  For a fixed $[a,m] = [a,p_{i}^{j} m']$ and $x \in Q_{i-1}$, there are at most $p_{i}^{\gamma_{i} - j}$ values of $y \in \mathcal A_{p_{i}^{\gamma_{i}}}$ such that $(x,y) \in [a,m]$ (as there are exactly $p_{i}^{\gamma_{i} - j}$ such values of $y \in \{ 0, 1, \ldots, p_{i}^{\gamma_{i}}-1 \}$).  From the definition of $\rho(p_{i}^{\gamma_{i}})$, we deduce that
\[
\alpha_{i}(x) 
\le  \dfrac{p_{i}^{\gamma_{i}}}{|\mathcal A_{p_{i}^{\gamma_{i}}}|} \sum_{j=1}^{\infty} \dfrac{1}{p_{i}^{j}} \sum_{\substack{[a, m] \in \mathcal C' \\ p_{i}^{j} \exactmid m, \,P(m) = p_{i} \\ x \in [a,m/p_{i}^{j}]}} 1
= \dfrac{1}{\rho(p_{i}^{\gamma_{i}})} \sum_{j=1}^{\infty} \dfrac{1}{p_{i}^{j}} \sum_{\substack{[a, p_{i}^{j} m'] \in \mathcal C' \\ P(m') < p_{i} \\ x \in [a,m']}} 1.
\]
Therefore, we obtain
\begin{equation}\label{expval2ndmomeq1}
\begin{split}
\mathbb{E}_{i-1}[\alpha_i^2]
&= \sum_{x \in Q_{i-1}} \mathbb{P}_{i-1}(x) \alpha_{i}^{2}(x) \\[5pt]
&\le \dfrac{1}{\rho(p_{i}^{\gamma_{i}})^{2}} \,\sum_{j_{1}=1}^{\infty} \,\sum_{j_{2}=1}^{\infty} \dfrac{1}{p_{i}^{j_{1}+j_{2}}} 
\sum_{\substack{[a_{1}, p_{i}^{j_{1}} m'_{1}] \in \mathcal C' \\ P(m'_{1}) < p_{i}}}
\,\sum_{\substack{[a_{2}, p_{i}^{j_{2}} m'_{2}] \in \mathcal C' \\ P(m'_{2}) < p_{i}}}
\,\sum_{\substack{x \in Q_{i-1} \\ x \in [a_{1},m'_{1}] \cap [a_{2},m'_{2}]}}  \mathbb{P}_{i-1}(x).
\end{split}
\end{equation}
This last sum is simply the probability (or weight) of $[a_{1},m'_{1}] \cap [a_{2},m'_{2}] \cap Q_{i-1}$, that is of the set of integers in $Q_{i-1} \subseteq [1,\ell_{i-1}]$ that are congruent to $a_{1} {\hskip -3pt}\pmod{m'_{1}}$ and congruent to $a_{2} {\hskip -3pt}\pmod{m'_{2}}$. 
If there are integers that are congruent to $a_{1} {\hskip -3pt}\pmod{m'_{1}}$ and congruent to $a_{2} {\hskip -3pt}\pmod{m'_{2}}$, then they correspond to a single residue class, say $a$, modulo $\lcm(m'_{1},m'_{2})$.
From \eqref{probdistPi}, we see that 
\begin{equation}
\label{Pixybound}
\mathbb P_{i}(x,y) \le \dfrac{\mathbb P_{i-1}(x)}{(1-\delta_{i}) |\mathcal A_{p_{i}^{\gamma_{i}}}|}
\qquad \text{for } 1 \le i \le r.
\end{equation}
If we write
\begin{equation}\label{lcmeq}
\lcm(m'_{1},m'_{2}) = \prod_{j=1}^{i-1} p_{j}^{e_{j}},
\end{equation}
where the $e_{j}$ are nonnegative integers, then we can view the $x \in [a_{1},m'_{1}] \cap [a_{2},m'_{2}] \cap Q_{i-1}$ as corresponding to the
\[
(x_{1},\ldots,x_{i-1}) \in \mathcal A_{p_{1}^{\gamma_{1}}} \times \mathcal A_{p_{2}^{\gamma_{2}}} \times \cdots \times \mathcal A_{p_{i-1}^{\gamma_{i-1}}},
\qquad \text{where $x_{j} \equiv a {\hskip -8pt}\pmod{p_{j}^{e_{j}}}$ for $1 \le j \le i-1$}.
\]
Note that if $e_{j} = 0$ for some $j$, then $x_{j}$ ranges over all values in $\mathcal A_{p_{j}^{\gamma_{j}}}$.
Taking into account \eqref{PQiQiminus1eq} in the case that $e_{j} = 0$ and \eqref{Pixybound} in the case $e_{j} > 0$, we obtain
\begin{align*}
\sum_{\substack{x \in Q_{i-1} \\ x \in [a_{1},m'_{1}] \cap [a_{2},m'_{2}]}}  \mathbb{P}_{i-1}(x)
&= \sum_{\substack{x_{1} \in \mathcal A_{p_{1}^{\gamma_{1}}} \\ x_{1} \in [a,p_{1}^{e_{1}}]}}
\,\sum_{\substack{x_{2} \in \mathcal A_{p_{2}^{\gamma_{2}}} \\ x_{2} \in [a,p_{2}^{e_{2}}]}}
\cdots \,\sum_{\substack{x_{i-1} \in \mathcal A_{p_{i-1}^{\gamma_{i-1}}} \\ x_{i-1} \in [a,p_{i-1}^{e_{i-1}}]}} \mathbb P_{i-1}(x_{1},\ldots,x_{i-1}) \\[5pt]
&\le \prod_{\substack{p_{j} \,\mid\, \lcm(m'_{1},m'_{2})\\ p_{j}^{e_{j}} \,\exactmid\, \lcm(m'_{1},m'_{2})}} \dfrac{p_{j}^{\gamma_{j}}}{p_{j}^{e_{j}}(1-\delta_{j}) |\mathcal A_{p_{j}^{\gamma_{j}}}|}
= \prod_{\substack{p_{j} \,\mid\, \lcm(m'_{1},m'_{2})\\ p_{j}^{e_{j}} \,\exactmid\, \lcm(m'_{1},m'_{2})}} \dfrac{1}{p_{j}^{e_{j}}(1-\delta_{j}) \rho(p_{j}^{\gamma_{j}})}.
\end{align*}
Next, we want to sum this product over the third and fourth sums on the right side of \eqref{expval2ndmomeq1}, but we do so by summing instead over the different possibilities for $m'_{1}$ and $m'_{2}$ where \eqref{lcmeq} holds for some fixed choice of nonnegative integers $e_{j}$.  Given \eqref{lcmeq}, we see that for each $p_{j}$, at least one of $p_{j}^{e_{j}} \exactmid m'_{1}$ and $p_{j}^{e_{j}} \exactmid m'_{2}$ holds and the other is divisible by $p_{j}^{t}$ for some $t \in \{ 0, 1, \ldots, e_{j} \}$.  As we only want to count the possibility that both $p_{j}^{e_{j}} \exactmid m'_{1}$ and $p_{j}^{e_{j}} \exactmid m'_{2}$ hold once, we see that there are a total of $(2e_{1}+1)(2e_{2}+1) \cdots (2e_{i-1}+1)$ possibilities for $m'_{1}$ and $m'_{2}$ satisfying \eqref{lcmeq}.  
Note that if $e_{j} = 0$, the factor $2e_{j}+1 = 1$ can be omitted.  
Setting $d = \lcm(m'_{1},m'_{2})$, we deduce then that
\begin{align*}
\sum_{\substack{[a_{1}, p_{i}^{j_{1}} m'_{1}] \in \mathcal C' \\ P(m'_{1}) < p_{i}}}
\,\sum_{\substack{[a_{2}, p_{i}^{j_{2}} m'_{2}] \in \mathcal C' \\ P(m'_{2}) < p_{i}}}
\,\sum_{\substack{x \in Q_{i-1} \\ x \in [a_{1},m'_{1}] \cap [a_{2},m'_{2}]}}  \mathbb{P}_{i-1}(x)
&\le \sum_{d \mid \ell_{i-1}} \,\prod_{\substack{p_{j}^{e_{j}} \exactmid d \\ e_{j} \ge 1}} \,\dfrac{2e_{j}+1}{p_{j}^{e_{j}}(1-\delta_{j}) \rho(p_{j}^{\gamma_{j}})}  \\[5pt]
&= \prod_{j < i} \bigg( 1 + \sum_{e_{j} = 1}^{\gamma_{j}} \,\dfrac{2e_{j}+1}{p_{j}^{e_{j}}(1-\delta_{j}) \rho(p_{j}^{\gamma_{j}})} \bigg) \\[5pt]
&\le \prod_{j < i} \bigg( 1 + \sum_{e_{j} = 1}^{\infty} \,\dfrac{2e_{j}+1}{p_{j}^{e_{j}}(1-\delta_{j}) \rho(p_{j}^{\gamma_{j}})} \bigg) \\[5pt]
&= \prod_{j < i} \bigg( 1 + \dfrac{3p_{j}-1}{(1-\delta_{j}) (p_{j}-1)^{2} \rho(p_{j}^{\gamma_{j}})} \bigg).
\end{align*}
Summing over $j_{1}$ and $j_{2}$ in \eqref{expval2ndmomeq1}, we obtain now that
\[
\mathbb{E}_{i-1}[\alpha_i^2] \le \dfrac{1}{\rho(p_{i}^{\gamma_{i}})^{2} (p_{i}-1)^{2}} \prod_{j < i} \bigg( 1 + \dfrac{3p_{j}-1}{(1-\delta_{j}) (p_{j}-1)^{2} \rho(p_{j}^{\gamma_{j}})} \bigg).
\]
From \eqref{measureBi}, we see that, for $\delta_i > 0$, we have
\begin{equation}\label{PiBi2ndmombd}
\mathbb{P}_i(B_i) \le \dfrac{1}{4\,\rho(p_{i}^{\gamma_{i}})^{2} \,\delta_i (1-\delta_i)\,(p_{i}-1)^{2}} \,\prod_{j < i} \bigg( 1 + \dfrac{3p_{j}-1}{(1-\delta_{j}) (p_{j}-1)^{2} \rho(p_{j}^{\gamma_{j}})} \bigg).
\end{equation}
With the above background, we are now ready to prove Theorem~\ref{densitythm}.

\begin{proof}[Proof of Theorem~\ref{densitythm}.]
We fix a set of congruences $\mathcal C'$ as above with the goal of showing that if the minimum modulus is large, there are integers in $[1,L_{r}] \cap \mathcal \mathcal S_{r}(\vec{v})$ which are not covered by congruences in $\mathcal C'$.
We set $\delta_{i} = 0$ for $i \le i_{0}$ where $i_{0}$ is to be determined but to be chosen large.  For $i > i_{0}$, we take $\delta_{i} = 1/2$.  
First, fix $i > i_{0}$.  The probability measure given by $\mathbb P = \mathbb P_{r}$ for the set of integers $n \in Q_{r}$ satisfying some congruence in $\mathcal C'$ with modulus $m$ satisfying $P(m) = p_{i}$ satisfies the bound given by \eqref{PiBi2ndmombd}.  
We obtain from Condition~A that there are constants $k_{j}$ depending at most on $c_{0}$ such that
\[
\prod_{j < i} \bigg( 1 + \dfrac{3p_{j}-1}{(1-\delta_{j}) (p_{j}-1)^{2} \rho(p_{j}^{\gamma_{j}})} \bigg)
< \prod_{j < i} \bigg( 1 + \dfrac{6p_{j}-2}{(p_{j}-1)^{2} c_{0}} \bigg) \le \prod_{j < i} \bigg( 1 + \dfrac{k_{1}}{p_{j}} \bigg)
\le \log^{k_{2}} p_{i} 
\]
so that \eqref{PiBi2ndmombd} implies
\[
\mathbb{P}_i(B_i) \le \dfrac{k_{3} \log^{k_{2}} p_{i}}{p_{i}^{2}}.
\]
We can now choose $i_{0}$ sufficiently large (depending only on $c_{0}$) such that
\begin{equation}\label{ultimate1}
    \sum_{i > i_{0}} \mathbb{P}_i(B_i) < \dfrac{1}{2}.
\end{equation}
We fix such an $i_{0}$.  
We may and do suppose $i_{0} \le r$; otherwise, $r < i_{0}$, a bound depending only on $c_{0}$, and an argument which replaces $i_{0}$ below with $r$ will complete the proof.

It suffices now to show that the set of integers in $Q_r$ covered by congruences in $\mathcal C'$ with modulus $m$ satisfying $P(m) \le p_{i_{0}}$ has probability measure given by $\mathbb P$ that is at most $1/2$.  We have not yet used that the moduli in $\mathcal C'$ are large, and this is where we do so.  
As noted earlier, since $\delta_{i} = 0$ for $i \le i_{0}$, for such $i$, the probability distribution $\mathbb P^{\uparrow}_{i}(x)$ is the uniform distribution on $Q_{r}$.  Suppose $[a,m] \in \mathcal C'$ with $P(m) = p_{i} \le p_{i_{0}}$.  Combining $\delta_{i} = 0$ with \eqref{probdistPi}, we obtain
\begin{align*}
\sum_{\substack{x \in Q_{r} \\ x \in [a,m]}} \mathbb P^{\uparrow}_{i}(x)
&= \sum_{\substack{x \in Q_{i} \\ x \in [a,m]}} \mathbb P_{i}(x)
= \sum_{\substack{x \in Q_{i} \\ x \in [a,m]}} \dfrac{1}{|\mathcal A_{p_{1}^{\gamma_{1}}}| \,|\mathcal A_{p_{2}^{\gamma_{2}}}| \cdots |\mathcal A_{p_{i}^{\gamma_{i}}}|} \\[5pt]
&\le \dfrac{\ell_{i}}{m} \cdot
\dfrac{1}{|\mathcal A_{p_{1}^{\gamma_{1}}}| \,|\mathcal A_{p_{2}^{\gamma_{2}}}| \cdots |\mathcal A_{p_{i}^{\gamma_{i}}}|}
\le \dfrac{\ell_{i_{0}}}{m} \cdot
\dfrac{1}{|\mathcal A_{p_{1}^{\gamma_{1}}}| \,|\mathcal A_{p_{2}^{\gamma_{2}}}| \cdots |\mathcal A_{p_{i_{0}}^{\gamma_{i_{0}}}}|}
= \dfrac{1}{m} \,\prod_{j = 1}^{i_{0}} \dfrac{1}{\rho(p_{j}^{\gamma_{j}})}.
\end{align*}
The last product is bounded above by $(1/c_{0})^{i_{0}}$, a bound depending only on $c_{0}$.
Since the sum of the reciprocals of the $p_{i_{0}}$-smooth numbers converges, there is a positive integer $k_{4}$ depending on $c_{0}$ such that \[
\bigg( \dfrac{1}{c_{0}} \bigg)^{i_{0}}
\sum_{\substack{m \ge k_{4} \\ m \text{ $p_{i_{0}}$-smooth}}} \dfrac{1}{m} < \dfrac{1}{2}.
\] 
Therefore, if all the congruences in $\mathcal C'$ have moduli at least as large as $k_{4}$, then the set of integers covered by congruences in $\mathcal C'$, with moduli $m$ satisfying $P(m) \le p_{i_{0}}$, has probability measure at most
\begin{equation}\label{ultimate2}
\sum_{i = 1}^{i_{0}} \mathbb P_{i}(B_{i})
\le \bigg( \dfrac{1}{c_{0}} \bigg)^{i_{0}}
\sum_{\substack{m \ge k_{4} \\ m \text{ $p_{i_{0}}$-smooth}}} \dfrac{1}{m} < \dfrac{1}{2}.
\end{equation}
Combining \eqref{ultimategoal}, \eqref{ultimate1}, \eqref{ultimate2} and Lemma~\ref{Srsuffices}, we see that some integers in $Q_{r}$ are not covered by congruences in $\mathcal C'$, finishing the proof of Theorem~\ref{densitythm} (with $K = k_{4}$).
\end{proof}

\bigskip\noindent
\textbf{Acknowledgment of Use of AI:}
After writing the paper, both ChatGPT-6.1 and the SAIR open math model were used to check for mathematical accuracy and grammar.  No mathematical errors were detected and only minor changes in the wording were implemented as a result.

\end{document}